\documentstyle[amssymb,preprint,aps]{revtex}

\tightenlines

\begin{document}
\draft
\author{D. Galetti, J.T. Lunardi\thanks{%
On leave from Departamento de Matem\'{a}tica e Estat\'{\i}{}stica, Setor de
Ci\^{e}ncias Exatas e Naturais, Universidade Estadual de Ponta Grossa, Ponta
Grossa -- PR -- Brazil.}, B.M. Pimentel and M. Ruzzi}
\address{Instituto de F\'{\i}{}sica Te\'{o}rica\\
Universidade Estadual Paulista - UNESP\\
Rua Pamplona 145\\
01405-900 S\~{a}o Paulo - SP - BRAZIL}
\title{Irreducibility and Compositeness in q-Deformed Harmonic Oscillator Algebras}
\maketitle

\begin{abstract}
q-Deformed harmonic oscillator algebra for real and root of unity values of
the deformation parameter is discussed by using an extension of the number
concept proposed by Gauss, namely the Q-numbers. A study of the reducibility
of the Fock space representation which explores the properties of the Gauss
polynomials is presented. When the deformation parameter is a root of unity,
an interesting result comes out in the form of a reducibility scheme for the
space representation which is based on the classification of the primitive
or non-primitive character of the deformation parameter. An application is
carried out for a $q$-deformed harmonic oscillator Hamiltonian, to which the
reducibility scheme is explicitly applied. For finite-dimensional spaces
associated to non-primitive roots of unity the compositeness of the
k-fermions/quons is discussed.
\end{abstract}

\pacs{Pacs: 03.65.Fd, 03.65.-w, 02.10.Lh\\
Keywords: q-deformed algebras, q-deformed harmonic oscillator, deformation
at roots of unity}

\section{Introduction}

In the last decades $q$-deformed algebras\cite
{drinfeld,jimbo,kulish,sklyanin} have been object of interest in the
literature and a great effort has been devoted to its understanding and
development\cite{chari,lohe,demichev}. In particular, the interest in
q-deformed algebras resides in the fact that they are deformed versions of
the standard Lie algebras, and give them back as the deformation parameter $%
q $ goes to unity. Furthermore, since it is known that the deformed algebras
encompass a set of symmetries that is richer than that of the Lie algebras,
one is tempted to recognize that quantum algebras can be the appropriate
tool to be dealt with in describing symmetries of physical systems which
cannot be properly treated within the Lie algebras, although the direct
interpretation of the deformation in these cases is sometimes incomplete or
even completely lacking. For instance, in some cases like the XXZ-model,
where the ferromagnetic/antiferromagnetic nature of a spin $\frac{1}{2}$
chain of length $N$ can be simulated through the introduction of a $q$%
-deformed algebra\cite{pasquier}, or the rotational bands in deformed nuclei
and molecules which can be fitted via a $q$-rotor Hamiltonian\cite
{iwao,bonatsos,celeghini}, instead of using the variable moment of inertia
(VMI model), the physical meaning of the deformation parameter is
established. Notwithstanding this interpretation difficulty, from the
original studies which appeared in connection with problems related to
solvable statistical mechanics models\cite{baxter} and quantum inverse
scattering theory\cite{fadeev}, a solid development has emerged which
encompass nowadays various branches of mathematical problems related to
physical applications, such as deformed superalgebras\cite{chaichian2}, knot
theories\cite{kauffman}, noncommutative geometries\cite{manin} and so on.
The introduction of a $q$-deformed bosonic harmonic oscillator is a subject
of great interest in this context and, as a tool for providing a boson
realisation of the quantum algebra $su_{q}(2)$, brought to light new
commutation relations\cite{biedenharn,macfarlane} which have been
extensively discussed in the literature.

On the other hand, some concepts directly related to the arithmetical
foundations of deformed algebras were well known to mathematicians since the
last century\cite{dickson,gauss,andrews}. For instance, the Gauss
polynomials appearing in restricted partition theory\cite{andrews} can be
directly interpreted as a $q$-generalization of the standard binomials; as
such, the Gauss polynomials, or the $q$-binomials, as they are sometimes
known, also generalize the concept of number as well. In that form, the
Gauss extension of the number concept, sometimes known as $Q$-number\cite
{andrews2,gasper}, is also related to the usual $q$-bracket of extensive use
in deformed algebras. In this connection, if, in general, the surprising
effectiveness of number theory seems not to be completely realized, the
success of recent examples pervading several areas can be credited to the
use of that branch of science: solvable models in statistical mechanics
benefited from Rogers-Ramanujan-Baxter relations, computation and
cryptography, the fourth test on general relativity, dynamical systems, and
primitive-roots-of-unity-based reflecting gratings in concert halls have
their very foundations on basic number theory and algorithms \cite{appl}.

In this paper we want to address the question whether the extension of the
number concept proposed by Gauss, namely the $Q$-numbers, can farther help
us in the study of the $q$-deformed harmonic oscillator. To this aim we are
directly guided by the central role played by the number concept in this
context. Based on this, we introduce the $Q$-numbers as our starting point
to define the action of the creation/annihilation operators on the Fock
space states. \ From this we show how we obtain a version of the $q$%
-deformed harmonic oscillator algebra already discussed in the literature 
\cite{arik,atak,kib}, the $A_{q}$ algebra. We also discuss how a second set
of operators obeying the $q$-deformed harmonic oscillator algebra can be
introduced, the $\overline{A_{q}}$ algebra, such that they satisfy the
conjugate relations with respect to the $A_{q}$ algebra, and discuss some
possible reductions of the algebra when we choose the allowed values of the
deformation parameter $q$. The cases for real and roots-of-unity values of $%
q $ are analysed. Furthermore, we also show how the reducibility of the
algebra space representation appears for the different values of $q.$

Using the algebras $A_{q}$ and $\overline{A_{q}}$ \ we introduce\ a
self-adjoint $q$-deformed harmonic oscillator Hamiltonian, akin to that
proposed by Floratos and Tomaras\cite{flor} and related to a system of two
anyons, which allows us to test the reducibility criteria discussed before.
This allows us to separate the physical systems according to the different
algebras obtained for the different values of the deformation parameter $q$.
In this form, we show how it is possible to distinguish different subsystems
within the original oscillator Hamiltonian when $q$ assume nonprimitive
roots of unity values. In this sense, we discuss the possibility of
uncovering the compositeness character of the so-called $k$-fermions when
discussing the reducibility of the representation space for $q$-deformed
oscillator algebra at the roots of unity.

This paper is organized as follows: Section II is devoted to a brief review
of the Gauss polynomials ($Q$-numbers) and their basic properties. In
section III we derive and discuss the $q$-oscillator algebras from the
extended number concept and in section IV we present the conditions for the
reducibility of the Fock space representation. $q$-Deformed oscillator
Hamiltonians are discussed in Section V, where examples of how the
reducibility conditions sieve the space representation into subspaces are
also exhibited. Finally the conclusions are presented in Section VI.

\section{Gauss Polynomials: $Q$-Numbers}

The generating function of restricted partitions of a positive integer $N$
into at most $m$ parts, each $\leq n,$\ is written as

\begin{equation}
G\left( n,m;q\right) =\frac{\left( 1-q^{n+m}\right) \left(
1-q^{n+m-1}\right) ...\left( 1-q^{m+1}\right) }{\left( 1-q\right) \left(
1-q^{2}\right) ...\left( 1-q^{n}\right) },
\end{equation}
$q\neq 1$, and the Gauss polynomials are defined through the relation 
\begin{equation}
\left[ 
\begin{tabular}{l}
$n$ \\ 
$m$%
\end{tabular}
\right] =G\left( n-m,m;q\right) ,
\end{equation}
which is valid for $0\leq m\leq n,$ and zero otherwise\cite{andrews}. The
Gauss polynomial is a polynomial of degree $m(n-m)$ in $q$ that presents a
very important property, namely 
\begin{equation}
\lim_{q\rightarrow 1}\left[ 
\begin{tabular}{l}
$n$ \\ 
$m$%
\end{tabular}
\right] =%
{n \choose m}%
,
\end{equation}
where $%
{n \choose m}%
$ is the standard binomial. Thus, we conclude that the Gauss polynomials
generalize the concept of binomials and, furthermore, as a special and
important case, with $m=1$, the Gauss polynomial, that is now denoted $Q$%
-number, extends the concept of number since 
\begin{equation}
\lim_{q\rightarrow 1}\left[ 
\begin{tabular}{l}
$n$ \\ 
$1$%
\end{tabular}
\right] =%
{n \choose 1}%
=n.
\end{equation}

On the other hand, this polynomial also allows us to establish inner contact
with some aspects of number theory, since when $q$ is a $n$th root of unity,
we have 
\begin{equation}
\left[ 
\begin{tabular}{l}
$n$ \\ 
$1$%
\end{tabular}
\right] =1+q+q^{2}+...+q^{n-1}=\frac{1-q^{n}}{1-q}=0.
\end{equation}
This is the fundamental equation whose $n$ solutions are roots of unity;
furthermore, for $n$ prime, $n-1$ of these are {\it primitive} roots \cite
{mathews}.

Besides those important properties, the Gauss polynomials also satisfy the
additional following relations\cite{andrews}

\begin{equation}
\left[ 
\begin{tabular}{l}
$n$ \\ 
$0$%
\end{tabular}
\right] =\left[ 
\begin{tabular}{l}
$n$ \\ 
$n$%
\end{tabular}
\right] =1,
\end{equation}

\begin{equation}
\left[ 
\begin{tabular}{l}
$n$ \\ 
$m$%
\end{tabular}
\right] =\left[ 
\begin{tabular}{l}
$n$ \\ 
$n-m$%
\end{tabular}
\right] ,
\end{equation}
\begin{equation}
\left[ 
\begin{tabular}{l}
$n$ \\ 
$m$%
\end{tabular}
\right] =\left[ 
\begin{tabular}{l}
$n-1$ \\ 
$m$%
\end{tabular}
\right] +q^{n-m}\left[ 
\begin{tabular}{l}
$n-1$ \\ 
$m-1$%
\end{tabular}
\right] ,
\end{equation}
\begin{equation}
\left[ 
\begin{tabular}{l}
$n$ \\ 
$m$%
\end{tabular}
\right] =\left[ 
\begin{tabular}{l}
$n-1$ \\ 
$m-1$%
\end{tabular}
\right] +q^{m}\left[ 
\begin{tabular}{l}
$n-1$ \\ 
$m$%
\end{tabular}
\right] .  \label{m}
\end{equation}

\section{$q$-Oscillator Algebras}

Let us consider, as our starting point, the standard Fock space generated by 
$\left\{ |n\rangle \right\} ,$%
\begin{equation}
a\mid 0\rangle =0,\;\;
\end{equation}
\begin{equation}
a^{\dagger }\mid n\rangle =\sqrt{n+1}\mid n+1\rangle \;,\;\;\;\;\;a\mid
n\rangle =\sqrt{n}\mid n-1\rangle ,
\end{equation}
and 
\begin{equation}
\hat{N}\mid n\rangle =n\mid n\rangle ,
\end{equation}
where the creation and annihilation operators obey the following commutation
relations 
\begin{equation}
aa^{\dagger }-a^{\dagger }a=1;\;\;\;\;\left[ \hat{N},a^{\dagger }\right]
=a^{\dagger };\;\;\;\;\;\left[ \hat{N},a\right] =-a\;,  \label{B}
\end{equation}
from which it follows that 
\begin{equation}
\hat{N}=a^{\dagger }a.
\end{equation}

Since the number concept is inherent to the Fock description, we are
strongly motivated by the results of the previous section to construct a new
pair of creation and annihilation operators in such a form to deal with that
generalized number concept. To this aim we introduce new operators, whose
matrix elements in the Fock space involve the Gauss polynomials 
\begin{equation}
a_{-}\mid 0\rangle =0
\end{equation}
\begin{equation}
a_{+}\mid n\rangle =\sqrt{\left\{ n+1\right\} _{q}}\mid n+1\rangle
\label{a+}
\end{equation}
\begin{equation}
a_{-}\mid n\rangle =\sqrt{\left\{ n\right\} _{q}}\mid n-1\rangle  \label{a-}
\end{equation}
\begin{equation}
\hat{N}\mid n\rangle =n\mid n\rangle
\end{equation}
\begin{equation}
\left[ \hat{N},a_{+}\right] =a_{+}  \label{n1}
\end{equation}
\begin{equation}
\left[ \hat{N},a_{-}\right] =-a_{-},  \label{n2}
\end{equation}
although $\hat{N}$ $\neq $ $a_{+}a_{-}$. Here we have adopted the notation 
\begin{equation}
\left\{ n\right\} _{q}\equiv \left[ 
\begin{tabular}{l}
$n$ \\ 
$1$%
\end{tabular}
\right] .
\end{equation}

We can pose now the question: what is the algebra satisfied by $a_{+}$ and $%
a_{-}$? Since 
\begin{equation}
a_{-}a_{+}\mid n\rangle =\left\{ n+1\right\} _{q}\mid n\rangle ,
\end{equation}
\begin{equation}
a_{+}a_{-}\mid n\rangle =\left\{ n\right\} _{q}\mid n\rangle ,
\end{equation}
and considering relation (\ref{m}), we conclude that 
\begin{equation}
a_{-}a_{+}-qa_{+}a_{-}=1,  \label{Aq}
\end{equation}
that is a q-deformed commutation relation as already exhibited in the
literature\cite{arik,atak,kib}. \ Let us denote \ relations \ (\ref{n1}, \ref
{n2}, \ref{Aq}) by $A_{q}$ algebra. We can construct an $\bar{A}_{q}$
algebra out of the relations conjugated to those defining the $A_{q}$
algebra (\ref{n1}, \ref{n2}, \ref{Aq}): 
\begin{equation}
\left[ \hat{N},a_{+}^{\dagger }\right] =-a_{+}^{\dagger },  \label{n1b}
\end{equation}
\begin{equation}
\left[ \hat{N},a_{-}^{\dagger }\right] =+a_{-}^{\dagger },  \label{n2b}
\end{equation}
\begin{equation}
a_{+}^{\dagger }a_{-}^{\dagger }-q^{\ast }a_{-}^{\dagger }a_{+}^{\dagger }=1.
\label{Aqb}
\end{equation}
In principle, these operators act on the dual space (bra space) to the
considered Fock (ket) space. However, we can infer the action of these
operators onto the ket space just by using the orthonormality of the states $%
|n\rangle $. It yields 
\begin{equation}
a_{-}^{\dagger }\mid n\rangle =\left( \sqrt{\left\{ n+1\right\} _{q}}\right)
^{\ast }\mid n+1\rangle ,  \label{a+d}
\end{equation}
and 
\begin{equation}
a_{+}^{\dagger }\mid n\rangle =\left( \sqrt{\left\{ n\right\} _{q}}\right)
^{\ast }\mid n-1\rangle .  \label{a-d}
\end{equation}
With the above results, it is possible to examine if there is an algebra
relating\ the creation/anihilation operators of the $A_{q}$ algebra and
their respective Hermitian conjugates, constituents of the $\bar{A}_{q}$
algebra. Using the action of these operators over the ket space, it is
possible to obtain the following relations: 
\begin{equation}
a_{-}a_{-}^{\dagger }=\left| \left\{ \hat{N}+1\right\} _{q}\right| ,
\end{equation}
\begin{equation}
a_{-}^{\dagger }a_{-}=\left| \left\{ \hat{N}\right\} _{q}\right| ,
\end{equation}
and similarly 
\begin{equation}
a_{+}^{\dagger }a_{+}=\left| \left\{ \hat{N}+1\right\} _{q}\right| ,
\label{bob1}
\end{equation}
\begin{equation}
a_{+}a_{+}^{\dagger }=\left| \left\{ \hat{N}\right\} _{q}\right| .
\label{bob2}
\end{equation}
Here we shall only consider cases when $q$ is real valued or a root of
unity, which are the most commonly found cases in the literature.

For real $q$ it is possible to verify that 
\begin{equation}
\left| \left\{ \hat{N}\right\} _{q}\right| =\left\{ \hat{N}\right\} _{q},
\end{equation}
which together with Eqs. (\ref{bob1}) and (\ref{bob2}), and the recurrence
relation of the Gauss polynomials, Eq. (\ref{m}), yields: 
\begin{equation}
a_{-}a_{-}^{\dagger }-qa_{-}^{\dagger }a_{-}=1.  \label{bob3}
\end{equation}
Similarly 
\begin{equation}
a_{+}^{\dagger }a_{+}-qa_{+}a_{+}^{\dagger }=1.
\end{equation}
In this case (real $q$), through Eqs. (\ref{a+}, \ref{a-d}), it is possible
to identify $a_{-}^{\dagger }\equiv a_{+}$.

When $q$ is the fundamental root of unity it can be, by its turn, verified
that 
\begin{equation}
\left| \left\{ \hat{N}\right\} _{q}\right| =\left[ \hat{N}\right] _{q^{1/2}},
\label{bob4}
\end{equation}
where 
\begin{equation}
\left[ X\right] _{q}=\frac{q^{X}-q^{-X}}{q-q^{-1}}  \label{brac}
\end{equation}
defines the $q$-bracket of $X$. Equation (\ref{bob4}), together with Eqs. (%
\ref{bob1}) and (\ref{bob2}) yields 
\begin{equation}
a_{-}a_{-}^{\dagger }-q^{\frac{1}{2}}a_{-}^{\dagger }a_{-}=q^{-\frac{\hat{N}%
}{2}}.  \label{bob5}
\end{equation}
Similarly 
\begin{equation}
a_{+}^{\dagger }a_{+}-q^{\frac{1}{2}}a_{+}a_{+}^{\dagger }=q^{-\frac{\hat{N}%
}{2}}.  \label{bob5b}
\end{equation}
The last two equations characterize the $q$-oscillator algebra introduced by
Biedenharn and McFarlane \cite{biedenharn,macfarlane}.

On the other hand, when $q$ is a root of unity, except the fundamental one,
Eq. (\ref{bob4}) is no longer valid, instead 
\begin{equation}
\left| \left\{ \hat{N}\right\} _{q}\right| =\left| \left[ \hat{N}\right]
_{q^{1/2}}\right| .
\end{equation}

Using the definition of the $q$-bracket, Eq. (\ref{brac}), when $q$ is a
general root of unity, $q_{j}=\exp (\frac{2\pi i}{m}j)$, a relation between $%
\left[ k\right] _{q_{j}^{\frac{1}{2}}}$ and its $m-$complementar $\left[ m-k%
\right] _{q_{j}^{\frac{1}{2}}}$, can be directly obtained 
\[
\left[ m-k\right] _{q_{j}^{\frac{1}{2}}}=\frac{e^{i\frac{\pi }{m}%
j(m-k)}-e^{-i\frac{\pi }{m}j(m-k)}}{e^{i\frac{\pi }{m}j}-e^{-i\frac{\pi }{m}%
j}}=\left( -1\right) ^{j-1}\frac{e^{i\frac{\pi }{m}jk}-e^{-i\frac{\pi }{m}jk}%
}{e^{i\frac{\pi }{m}j}-e^{-i\frac{\pi }{m}j}} 
\]
\begin{equation}
\left[ m-k\right] _{q_{j}^{\frac{1}{2}}}=\left( -1\right) ^{j-1}\left[ k%
\right] _{q_{j}^{\frac{1}{2}}}.  \label{rel1}
\end{equation}
When $q$ is the fundamental root of unity then $j=1$, and we have 
\begin{equation}
\left[ m-k\right] _{q_{1}^{\frac{1}{2}}}=\left[ k\right] _{q_{1}^{\frac{1}{2}%
}}.  \label{rel2}
\end{equation}
Furthermore, for the case of the inverse of such root of unity, $%
q_{j}^{-1}=\exp \left( -\frac{2\pi i}{m}j\right) =\exp \left[ \frac{2\pi i}{m%
}(m-j)\right] $, we can verify in exactly the same way that 
\begin{equation}
\left[ k\right] _{q_{j}^{-\frac{1}{2}}}=\left( -1\right) ^{k-1}\left[ k%
\right] _{q_{j}^{\frac{1}{2}}}.  \label{rel3}
\end{equation}
Now, using Eqs. (\ref{rel1}-\ref{rel3}) we obtain the following additional
relation 
\begin{equation}
\left[ m-k\right] _{q_{j}^{-\frac{1}{2}}}=\left( -1\right) ^{m-k-1}\left[ m-k%
\right] _{q_{j}^{\frac{1}{2}}}.  \label{rel5}
\end{equation}
These relations will be shown to be useful when we deal with $q$-oscillator
Hamiltonians in finite-dimensional spaces.

\section{Reducibility of the Fock Representation}

Now, considering the actions of $a_{+}$ and $a_{-}$ on the Fock
representation, Eqs. (\ref{a+}) and (\ref{a-}), we will analyse its
reducibility properties. The various possibilities are studied below.

\subsection{First case: $\left\{ n\right\} _{q}\neq 0,\ \forall \ n>0$}

All states of the $\{|n\rangle \}$ representation can be obtained through
successive applications of $a_{+}$ over the vacuum. In that case, $%
\{|n\rangle \}$ is irreducible with respect to the algebra $\{a_{-},a_{+},%
\hat{N},I\}.$

\subsection{Second case: $\left\{ m\right\} _{q}=0,\ \left\{ n\right\}
_{q}\neq 0\ \;\forall \ n,\ 0<n<m$}

In that case 
\begin{equation}
a_{+}|m-1\rangle =0,  \label{13}
\end{equation}
and also 
\begin{equation}
a_{-}|m\rangle =0.  \label{14}
\end{equation}
From these results it follows that the subspace generated by $\{|0\rangle
,|1\rangle ,...,|m-1\rangle \}$ is invariant under the action of the set $%
\{a_{-},a_{+},\hat{N},I\}$, and it is then an {\it irrep} of dimension $m$
of the deformed algebra. For all $q\neq 1$,{\it \ i.e.,} deformed cases, the
hypothesis $\left\{ m\right\} _{q}=0,\ $ $\left\{ n\right\} _{q}\neq 0,\
\forall \ n,\ 0<n<m$\ can be written as 
\begin{equation}
\frac{q^{m}-1}{q-1}=0,\qquad \frac{q^{n}-1}{q-1}\neq 0,\qquad \forall \ n,\
0<n<m  \label{15}
\end{equation}

\begin{equation}
q^{m}=1,\qquad q^{n}\neq 1,\qquad \forall \ n,\ 0<n<m,  \label{16}
\end{equation}
which is the definition of the primitive $m$th roots of unity. Therefore
there will always be {\it irreps} of dimension $m$ whenever $q$ is a
primitive $m$th root of unity.

\subsection{Third case: $\exists \ l,\ 0<l<m\ \diagup \ \left\{ m\right\}
_{q}=0,\ \left\{ l\right\} _{q}=0\ \ $}

This is the equivalent to state that $q$ is a non-primitive root of unity.
Let us also suppose that $l$ is the smallest integer satisfying the
hypothesis above, {\it i.e., }it{\it \ }is the smallest number for which $%
q^{l}=1.$ Then 
\begin{equation}
q^{k}\neq 1,\qquad \forall \ \ k,\ 0<k<l,  \label{17}
\end{equation}
and therefore the subspace generated by $\{|0\rangle ,|1\rangle
,...,|m-1\rangle \}$ is reducible in {\it irreps} of dimension $l.$

Labelling the $m-1$ roots of unity as 
\begin{equation}
q_{j}=e^{2\pi i\frac{j}{m}},\qquad j=1,2,...,m-1,  \label{18}
\end{equation}
and for $r$ the greatest common divisor (GCD) of $m$ and $j$, {\it i.e.,} 
\begin{eqnarray}
j &=&rs  \nonumber \\
m &=&rl,\qquad  \label{19}
\end{eqnarray}
where $s/l$ is an irreducible fraction, then 
\begin{equation}
q_{j}=e^{2\pi i\frac{s}{l}}.  \label{20}
\end{equation}
In this way, the subspace generated by $\{|0\rangle ,|1\rangle
,...,|m-1\rangle \}$ is reducible, as we saw, in {\it irreps} of dimension $%
l $, which is the smallest value for which $q_{j}^{l}=1$.

Then, for each $m$th root of unity labelled by $j$, the dimension of the 
{\it irreps} will be $l=m/r$, where $r$ is the GCD of $j$ and $m$, and the $%
m-$dimensional representation breaks into $r$ {\it irreps} of dimension $%
l=m/r$.

\section{$q$-Deformed Oscillator Hamiltonian}

We can obtain a $q$-deformed Hermitian oscillator Hamiltonian from the
deformed operator algebra presented above through the direct construction 
\begin{equation}
H=\frac{1}{2}\hbar \omega \left( a_{-}a_{-}^{\dagger }+a_{-}^{\dagger
}a_{-}\right) =\frac{1}{2}\hbar \omega \left( a_{+}a_{+}^{\dagger
}+a_{+}^{\dagger }a_{+}\right) ,  \label{ham}
\end{equation}
that can be written, in general, as 
\begin{equation}
H=\frac{1}{2}\hbar \omega \left( \sqrt{\left\{ \hat{N}+1\right\} _{q}\left\{ 
\hat{N}+1\right\} _{q^{\ast }}}+\sqrt{\left\{ \hat{N}\right\} _{q}\left\{ 
\hat{N}\right\} _{q^{\ast }}}\right) ,
\end{equation}
which is equivalent to 
\begin{equation}
H=\frac{1}{2}\hbar \omega \left( \left| \left\{ \hat{N}+1\right\}
_{q}\right| +\left| \left\{ \hat{N}\right\} _{q}\right| \right) .
\end{equation}
As was discussed in the preceding sections 
\begin{equation}
\left| \left\{ \hat{N}\right\} _{q}\right| =\left\{ 
\begin{tabular}{l}
$\left\{ \hat{N}\right\} _{q},\qquad $for $q$ a real number \\ 
$\left| \left[ \hat{N}\right] _{q^{1/2}}\right| ,\qquad $for $q$ a root of
unity.
\end{tabular}
\right.
\end{equation}
So, for real $q$, the Hamiltonian, Eq. (\ref{ham}), is written as 
\begin{equation}
H=\frac{1}{2}\hbar \omega \left( \left\{ \hat{N}+1\right\} _{q}+\left\{ \hat{%
N}\right\} _{q}\right) ,
\end{equation}
and, for $q$ being a root of unity, it can be easily seen to reduce to 
\begin{equation}
H=\frac{1}{2}\hbar \omega \left( \left| \left[ \hat{N}+1\right]
_{q_{j}^{1/2}}\right| +\left| \left[ \hat{N}\right] _{q_{j}^{1/2}}\right|
\right) .  \label{qham}
\end{equation}
However, when $q$\ is furthermore singled out as the fundamental primitive
root of unity, the above expression, according to Eq. (\ref{bob4}), reduces
to 
\begin{equation}
H=\frac{1}{2}\hbar \omega \left( \left[ \hat{N}+1\right] _{q_{1}^{1/2}}+%
\left[ \hat{N}\right] _{q_{1}^{1/2}}\right) ,
\end{equation}
which is the usual proposal for the $q$-deformed oscillator \cite{biedenharn}%
. This last expression has the symmetry $q\rightarrow q^{-1}$, since this is
a symmetry of the bracket itself.{\bf \ }

Since the deformed oscillator Hamiltonian is written directly in terms of
the brackets of the operator $\hat{N}$, the reducibility properties of the
Fock representation space will appear in the spectrum of that operator as
well. In this sense, the spectrum will be broken into blocks associated to
subspaces of prime dimension whenever the initial space dimension is a
composite integer number and we work with the {\it non-primitive roots of
unity}.

As an application of what has been presented above we will discuss some
simple cases. In this connection, we need not to work with all the roots of
unity due to the properties presented for the brackets in the previous
sections. In fact, using relations (\ref{rel1}-\ref{rel5}), we need not
calculate the matrices representing the Hamiltonian for some primitive roots.

First, let us consider $m=2$. In this case, the matrix representing the $q$%
-deformed oscillator Hamiltonian is directly written since we only have to
work with the fundamental primitive root of unity, $q_{1}^{1/2}=\exp \left( i%
\frac{\pi }{2}\right) $. Using the fact that $\left[ 2\right]
_{q_{1}^{1/2}}=0$, we get 
\begin{equation}
H=\frac{1}{2}\hbar \omega \left( 
\begin{array}{cc}
\left[ 1\right] _{q_{1}^{1/2}} & 0 \\ 
0 & \left[ 1\right] _{q_{1}^{1/2}}
\end{array}
\right) =\frac{1}{2}\hbar \omega \left( 
\begin{array}{cc}
1 & 0 \\ 
0 & 1
\end{array}
\right) .
\end{equation}
For $m=3$, which is the next prime number, and also using Eq. (\ref{rel2}),
we get for the fundamental primitive root of unity, $q_{1}^{1/2}=\exp \left(
i\frac{\pi }{3}\right) ,$%
\begin{equation}
H=\frac{1}{2}\hbar \omega \left( 
\begin{array}{ccc}
\left[ 1\right] _{q_{1}^{1/2}} & 0 & 0 \\ 
0 & \left[ 1\right] _{_{q_{1}^{1/2}}}+\left[ 2\right] _{q_{1}^{1/2}} & 0 \\ 
0 & 0 & \left[ 2\right] _{q_{1}^{1/2}}
\end{array}
\right) =\frac{1}{2}\hbar \omega \left( 
\begin{array}{ccc}
1 & 0 & 0 \\ 
0 & 2 & 0 \\ 
0 & 0 & 1
\end{array}
\right) .
\end{equation}

If we now consider the case $m=6$, we can verify how the matrix representing
the Hamiltonian breaks into blocks,{\it \ }each with a prime dimension, as
occurs in the representation space of the $q$-deformed algebra. To this end,
let us first of all consider the Hamiltonian associated to the fundamental
primitive root of unity, $q_{1}^{1/2}=\exp \left( i\frac{\pi }{6}\right) $.
In this case, using Eq. (\ref{rel1}), we see that the matrix is also
symmetric and has the form 
\begin{equation}
H=\frac{1}{2}\hbar \omega \left( 
\begin{array}{cccccc}
1 &  &  &  &  &  \\ 
& 1+\left[ 2\right] _{q_{1}^{1/2}} &  &  & \text{{\Huge 0}} &  \\ 
&  & \left[ 2\right] _{q_{1}^{1/2}}+\left[ 3\right] _{q_{1}^{1/2}} &  &  & 
\\ 
&  &  & \left[ 3\right] _{q_{1}^{1/2}}+\left[ 2\right] _{q_{1}^{1/2}} &  & 
\\ 
& \text{{\Huge 0}} &  &  & \left[ 2\right] _{q_{1}^{1/2}}+1 &  \\ 
&  &  &  &  & 1
\end{array}
\right) .
\end{equation}

Now, if we consider the non-primitive roots of $m=6$, we see that there are
three of them, namely, $q_{2}^{1/2}=\exp \left( i\frac{\pi }{6}2\right) $, $%
q_{3}^{1/2}=\exp \left( i\frac{\pi }{6}3\right) $ and $q_{4}^{1/2}=\exp
\left( i\frac{\pi }{6}4\right) $ respectively. In fact, $q_{2}$ is the
inverse of $q_{4}$ and $q_{3}$ is its own inverse. For the first root the
Hamiltonian matrix will be 
\begin{equation}
H=\frac{1}{2}\hbar \omega \left( 
\begin{array}{cccccc}
1 &  &  &  &  &  \\ 
& 2 &  &  & \text{{\Huge 0}} &  \\ 
&  & 1 &  &  &  \\ 
&  &  & 1 &  &  \\ 
& \text{{\Huge 0}} &  &  & 2 &  \\ 
&  &  &  &  & 1
\end{array}
\right) ,
\end{equation}
which breaks into two blocks, each one being the matrix associated to a $m=3$
$q$-deformed oscillator. On the other hand, for the second non-primitive
root of unity, $q_{3}^{1/2}=\exp \left( i\frac{\pi }{6}3\right) $, we get 
\begin{equation}
H=\frac{1}{2}\hbar \omega \left( 
\begin{array}{cccccc}
1 &  &  &  &  &  \\ 
& 1 &  &  & \text{{\Huge 0}} &  \\ 
&  & 1 &  &  &  \\ 
&  &  & 1 &  &  \\ 
& \text{{\Huge 0}} &  &  & 1 &  \\ 
&  &  &  &  & 1
\end{array}
\right) .
\end{equation}
The three blocks associated to the $m=2\;q$-deformed oscillator are readily
seen in this case. Since the Hamiltonian is given by (\ref{qham}), and using
Eq. (\ref{rel3}), we can conclude that the same matrices would be obtained
if the inverse roots of unity were used. Therefore, for prime dimension
spaces the matrices representing the deformed oscillator Hamiltonian are
irreducible for any primitive root of unity. For nonprime integer space
dimension and deformations at the {\it non-primitive roots of unity}, the $q$%
-deformed oscillator represents in fact a composite system with as many
irreducible constituents (diagonal blocks) as are the number of prime
factors of the starting space dimension.

\section{Conclusions}

In the present paper, starting from the Gauss extension of the number
concept, we have reobtained the $q$-deformed harmonic oscillator algebra
discussed in \cite{arik,atak} for general deformation parameter $q$. For the
particular case of $q$ being a fundamental root of unity, we recover the
deformed harmonic oscillator algebra satisfied by $a_{-(+)}$ and $%
a_{-(+)}^{\dagger }$ as introduced by Biedenharn and MacFarlane\cite
{biedenharn,macfarlane}. On the other hand, some useful relations between
the Gauss polynomials and the standard $q$-bracket have also been discussed
for $q$ a root of unity.

A discussion on the dimensions of the Fock state space representation for $\
q$ real or a root of unity shows that they can be infinite as well as
finite-dimensional depending on $q$ being real or a root of unity, as it has
been already pointed out. However, as a further conclusion, it is also shown
that, for the particular case of $q$ being selected as a non-primitive root
of unity, the representation space, besides being finite-dimensional, also
breaks into subspaces, the dimension of each block being clearly defined by
the prime decomposition of the number characterizing the original space
dimension. In this form, it is shown that the use of non-primitive roots of
unity allows one to verify the reducibility character of the Fock space
representation, which, by its turn, shows that the subspaces characterized
by prime dimensions play the role of fundamental blocks within the full
space.

A $q$-deformed harmonic oscillator Hamiltonian was presented which allowed
us to fully exploit the Fock space reducibility discussed previously. For
the cases when $\ q$ was a primitive root of unity, the Hamiltonian matrix
only exhibited the usual symmetries of the Q-numbers, while for $q$ being a
non-primitive root of unity (which will occur only when the space dimension
is a composite number) the matrices reduced to submatrices along the
diagonal, thus indicating that the original $q$-deformed oscillator is in
fact made up of irreducible subsystems, each one of them of a prime
dimension. This result strongly suggests the conclusion that the so called
k-fermions, or quons, discussed in the context of roots of unity
deformation, are not necessarily fundamental entities, but they may be, in
some cases, composite systems made up of entities of more fundamental
character. This characterization can be directly verified by studying the
degree of reducibility of the space representation through the prime
decomposition of the space dimension from which we started.

{\bf Acknowledgement}: D.G. and B.M.P are partially supported by Conselho
Nacional de Desenvolvimento Cient\'{\i}{}fico e Tecnol\'{o}gico, CNPq;
J.T.L. is partially supported by PICDT/CAPES, and M.R. has a fellowship from
Funda\c{c}\~{a}o de Amparo \`{a} Pesquisa do Estado de S\~{a}o Paulo, FAPESP.

\appendix

\section{Polychronakos Realization}

We start from the fundamental relation (\ref{Aq}) which clearly reduces to
the classical oscillator algebra when $q\rightarrow 1$. We recall the
Polychronakos realization \cite{poly}: 
\begin{eqnarray}
a_{-} &=&U_{-}(q,\hat{N})a  \nonumber \\
a_{+} &=&U_{+}(q,\hat{N})a^{\dagger },  \label{clare}
\end{eqnarray}
where $\hat{N}$ is the usual nondeformed number operator. Using Eq. (\ref
{clare}) in Eq. (\ref{Aq}) we obtain: 
\begin{equation}
F(q,\hat{N}+1)-qF(q,\hat{N})=1,  \label{3}
\end{equation}
where 
\begin{equation}
F(q,\hat{N})=U_{+}(q,\hat{N})U_{-}(q,\hat{N}-1)\hat{N}.  \label{4}
\end{equation}
Representing Eq. (\ref{4}) on the Fock space $\{|n\rangle \}$ we get: 
\begin{equation}
F(q,n+1)-qF(q,n)=1,
\end{equation}
which is the recurrence relation (\ref{m}) for the Gauss polynomials. We
then may infer 
\begin{eqnarray}
a_{+}a_{-} &=&g(q,\hat{N})  \label{6.5} \\
a_{-}a_{+} &=&g(q,\hat{N}+1),  \label{7}
\end{eqnarray}
or in terms of $F$%
\begin{equation}
F(q,\hat{N})=g(q,\hat{N}).
\end{equation}

In order to fulfil the deformed algebra, it is enough that 
\begin{equation}
U_{+}(q,\hat{N})U_{-}(q,\hat{N}-1)\hat{N}=g(q,\hat{N}).  \label{9}
\end{equation}
Choosing 
\begin{equation}
U_{-}(q,\hat{N}-1)=\sqrt{\frac{\left\{ \hat{N}\right\} _{q}}{\hat{N}}}%
,\qquad U_{-}(q,\hat{N})=\sqrt{\frac{\left\{ \hat{N}+1\right\} _{q}}{\hat{N}%
+1}},  \label{10}
\end{equation}
we obtain 
\begin{equation}
U_{+}(q,\hat{N})=\sqrt{\frac{\left\{ \hat{N}\right\} _{q}}{\hat{N}}}.
\label{11}
\end{equation}
Therefore 
\begin{mathletters}
\label{12}
\begin{eqnarray}
a_{-} &=&a\sqrt{\frac{\left\{ \hat{N}\right\} _{q}}{\hat{N}}}  \label{12a} \\
a_{+} &=&a^{\dagger }\sqrt{\frac{\left\{ \hat{N}+1\right\} _{q}}{\hat{N}+1}.}
\label{12b}
\end{eqnarray}
This choice guarantees the unitarity $(a_{+}=a_{-}^{\dagger })$ when $q$ is
a real parameter. Otherwise, the representation turns out to be non-unitary.

\end{mathletters}

\end{document}